\documentclass{ifacconf}

\usepackage{epsfig}
\usepackage{graphicx}
\usepackage{subfigure}
\usepackage{psfrag}

\usepackage{amsmath}
\usepackage{amssymb}
\usepackage{amsbsy}
\usepackage{subeqnarray}

\usepackage{textcomp}

\usepackage{natbib}

\newtheorem{theorem}{Theorem}

\newtheorem{assumption}{Assumption}

\newtheorem{remark}{Remark}

\newcommand{\mbf}[1]
  {{\mathbf #1}}
\newcommand{\mbb}[1]
  {{\mathbb #1}}

\newcommand{\mbg}[1]
  {{\mbox{\boldmath$#1$}}}

\newcommand{\low}[1]
  {_{\scriptscriptstyle #1}}

\begin{document}

\begin{frontmatter}

\title{A Sensor-based Long Baseline Position and Velocity Navigation Filter for Underwater Vehicles}

\author[ALL]{Pedro~Batista} 
\author[ALL]{Carlos~Silvestre} 
\author[ALL]{Paulo~Oliveira} 

\address[ALL]{Instituto Superior T\'ecnico\\Av. Rovisco Pais, 1049-001 Lisboa, Portugal\\{\small \{pbatista,cjs,pjcro@isr.ist.utl.pt\}@isr.ist.utl.pt}}

\begin{abstract}
This paper presents a novel Long Baseline (LBL) position and
velocity navigation filter for underwater vehicles based directly on
the sensor measurements. The solution departs from previous
approaches as the range measurements are explicitly embedded in the
filter design, therefore avoiding inversion algorithms. Moreover,
the nonlinear system dynamics are considered to their full extent
and no linearizations are carried out whatsoever. The filter error
dynamics are globally asymptotically stable (GAS) and it is shown,
under simulation environment, that the filter achieves similar
performance to the Extended Kalman Filter (EKF) and outperforms
linear position and velocity filters based on algebraic estimates of
the position obtained from the range measurements.
\end{abstract}


\end{frontmatter}

\section{Introduction}
\label{sec:Introduction}

Accurate navigation systems are essential for the successful
operation of autonomous vehicles. Although there exist alternatives
such as terrain-based navigation, most navigation systems contain
an Inertial Navigation System (INS) that provides the state of the
vehicle by integrating, in open-loop, the information provided by
inertial sensors, e.g., accelerometers and rate gyros. Although INS
provides very good short term results, its performance necessarily
degrades over time, not only due to the integration of sensor noise
but also due to sensor bias errors. In order to overcome these
drawbacks, aiding devices are considered to correct INS errors. This
paper addresses the problem of vehicle navigation using ranges to a
set of landmarks disposed in a Long Baseline (LBL) configuration.

Among the myriad of aiding devices, the Global Positioning System
(GPS) is a very popular choice, see, e.g.,
\cite{paper:Sukkarieh:INSGPSGround}
\cite{paper:Yun:TestingEvaluationGPS-INSSystemSmallAUVNavigation},
and
\cite{paper:Fossen:NonlinearObserverGPS-INSIntegration}. For
underwater vehicles, GPS is not a solution due to the strong
attenuation that the electromagnetic field suffers in water. Therefore, other solutions have been sought
in the past, including acoustic positioning systems like LBL and
Short Baseline (SBL), see 
\cite{paper:Jouffroy:UnderwaterVehicleTrajectoryEstimationUsingContractingPDE-basedObservers},
\cite{paper:Kinsey:PreliminaryDVLNAVIntegratedNavigationSystemMannedUnmannedSubmersibles},
\cite{thesis:Larsen},
\cite{paper:Larsen:SLBL},
\cite{paper:Vaganay:ComparisonFixComputationFilteringAutonomousAcousticNavigation},
and references therein.
In \cite{patent:Youngberg:MethodGPSUnderwater}
the author proposes a GPS-like system consisting of
buoys equipped with DGPS receptors. A related solution, denominated
as GPS Intelligent Buoy (GIB) system, is now commercially
available, see \cite{paper:Thomas:GIBBuoys}. Further work on the GIB
underwater positioning system can be found in
\cite{paper:Alcocer:EKFGIBUnderwaterPositioningSystem}.
Position and linear velocity globally asymptotically stable (GAS) filters based on a Ultra-Short
Baseline (USBL) positioning system were presented by the authors in
\cite{paper:Batista:IFACWC:2008:2}.
For interesting discussions and more detailed surveys on underwater 
vehicle navigation and sensing devices see
\cite{paper:Whitcomb:2006:SurveyUnderwaterNavigation}
and
\cite{report:Leonard:AutonomousUnderwaterVehicleNavigation}.

This paper presents a position and linear velocity navigation
filter based on range measurements disposed in a LBL
configuration. Traditional solutions resort either to the
Extended Kalman Filter (EKF) or to solutions based on position
algebraic estimates obtained from the range measurements. The
solution presented in the paper departs from previous approaches
as the range measurements are explicitly embedded in the filter
design, therefore avoiding inversion algorithms. Moreover, the
nonlinear system dynamics are considered to their full extent and no
linearizations are carried out whatsoever, which allows to show that
the filter error dynamics are globally asymptotically stable.
Central to the proposed filtering framework is the derivation of a
linear time-varying (LTV) system that captures the dynamics of the
nonlinear system. The LTV model is achieved through
appropriate state augmentation, which is shown to mimic the nonlinear
system. Applications of the proposed solution are many and, under
simulation environment, it is shown that the filter achieves similar
performance to the Extended Kalman Filter (EKF) and outperforms
linear position and velocity filters based on position algebraic
estimates obtained directly from the range measurements.

The paper is organized as follows. The problem statement and 
nominal system dynamics are introduced in Section \ref{sec:PS},
while the filter design is detailed in Section \ref{sec:FD}.
Simulation results are presented in Section \ref{sec:SR} and Section
\ref{sec:C} summarizes the main conclusions of the paper.

\subsection{Notation}
Throughout the paper the symbol $\mbf{0}_{n \times m}$ denotes an
$n \times m$ matrix of zeros, $\mbf{I}_n$ an identity matrix with
dimension $n \times n$, and
$\rm{diag}(\mbf{A}_1, \ldots, \mbf{A}_n)$ a block diagonal matrix.
When the dimensions are omitted the matrices are assumed of
appropriate dimensions. If $\mbf{x}$ and $\mbf{y}$ are two vectors
of identical dimensions, $\mbf{x} \times \mbf{y}$ and
$\mbf{x} \cdot \mbf{y}$ represent the cross and inner products,
respectively.

\section{Problem Statement}
\label{sec:PS}

\begin{figure}[htbp]
	\center
	{
		\includegraphics[width=0.35\textwidth, keepaspectratio]{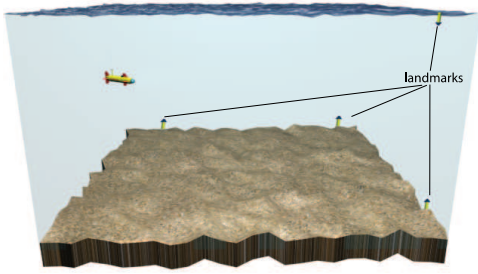}
	}
	\caption{Mission Scenario}
	\label{fig:PS:MS}
\end{figure}

Consider an underwater vehicle moving in a scenario where there is a
set of fixed landmarks installed in a Long Baseline configuration
and suppose that the vehicle measures the range to each of the
landmarks, as depicted in Fig. \ref{fig:PS:MS}.
Further assume that the vehicle is equipped with an Inertial
Measurement Unit (IMU), consisting of two triads of orthogonally
mounted accelerometers and rate gyros, and an Attitude and Heading
Reference System (AHRS). The problem considered in the paper is the
design of a sensor-based globally asymptotically stable filter to
estimate the position and linear velocity of the vehicle.

\subsection{System Dynamics}
\label{sec:PS:SD}

In order to detail the problem framework, let $\{I\}$ denote an
inertial reference coordinate frame and $\{B\}$ a coordinate frame
attached to the vehicle, commonly denominated as the body-fixed
coordinate frame. The linear motion of the vehicle is described by
\begin{equation}
	\dot{\mbf{p}}(t) = \mbf{R}(t) \mbf{v}(t),
	\label{eq:PS:SD:Der p}
\end{equation}
where $\mbf{p} \in \mbb{R}^3$ denotes the inertial position of the
vehicle, $\mbf{v} \in \mbb{R}^3$ is the velocity of the vehicle
relative to $\{I\}$ and expressed in body-fixed coordinates, and
$\mbf{R}$ is the rotation matrix from $\{B\}$ to $\{I\}$, which
satisfies
$\dot{\mbf{R}}(t) = \mbf{R}(t) \mbf{S}\left(\mbg{\omega}(t)\right)$,
where $\mbg{\omega} \in \mbb{R}^3$ is the angular velocity of
$\{B\}$, expressed in body-fixed coordinates, and
$\mbf{S}\left(\mbg{\omega}\right)$ is the skew-symmetric matrix such
that $\mbf{S}\left(\mbg{\omega}\right) \mbg{x}$ is the cross product
$\mbg{\omega} \times \mbg{x}$.

The AHRS provides the attitude rotation matrix $\mbf{R}$ and the
angular velocity $\mbg{\omega}$, while the IMU, assumed mounted at the
center of mass of the vehicle, measures the linear acceleration
$\mbf{a}$, which satisfies
\begin{equation}
	\mbf{a}(t) = \dot{\mbf{v}}(t) + \mbf{S}\left(\mbg{\omega}(t)\right) \mbf{v}(t) - \mbf{g}(t),
	\label{eq:PS:SD:a}
\end{equation}
where $\mbf{g} \in \mbb{R}^3$ denotes the acceleration of gravity,
expressed in body-fixed coordinates. Ideal accelerometers would not
measure the gravitational term but in practice this term must be
considered due to the inherent physics of the accelerometers, see
\cite{report:Kelly:NavigationSystems} for further details. Since the
magnitude of $\mbf{g}$ is usually well known, it would be possible
to cancel this term in \eqref{eq:PS:SD:a} as the attitude of the
vehicle is measured. However, even small errors on $\mbf{R}$ would
lead to large errors in the acceleration compensation. Therefore,
the acceleration of gravity is considered here as unknown state,
in addition to $\mbf{p}$ and $\mbf{v}$. The term
$\mbf{S}\left(\mbg{\omega}(t)\right) \mbf{v}(t)$ corresponds to the
Coriolis acceleration and cannot be neglected, particularly for
vehicles that execute aggressive maneuvers. Finally, let
$\mbf{s}_i \in \mbb{R}^3, \, i=1,\, \ldots,\, n\low{L},$ denote the
inertial positions of the landmarks. Then, the range measurements
are given by
\begin{equation}
	r_i(t) = \left\| \mbf{s}_i - \mbf{p}(t) \right\|.
	\label{eq:PS:SD:r_i}
\end{equation}

The derivative of the velocity can be written, from
\eqref{eq:PS:SD:a}, as
\begin{equation}
	\dot{\mbf{v}}(t) = \mbf{a}(t) - \mbf{S}\left(\mbg{\omega}(t)\right) \mbf{v}(t) + \mbf{g}(t).
	\label{eq:PS:SD:Der v}
\end{equation}
The derivative of $\mbf{g}$, assuming that the acceleration of
gravity is constant in inertial coordinates, is given by
\begin{equation}
	\dot{\mbf{g}}(t) = - \mbf{S}\left(\mbg{\omega}(t)\right) \mbf{g}(t).
	\label{eq:PS:SD:Der g}
\end{equation}
Combining \eqref{eq:PS:SD:Der p} and
\eqref{eq:PS:SD:r_i}-\eqref{eq:PS:SD:Der g} yields the nonlinear
system
\begin{equation}
	\left\{
		\begin{array}{l}
			\dot{\mbf{p}}(t) = \mbf{R}(t) \mbf{v}(t)\\
			\dot{\mbf{v}}(t) = \mbf{a}(t) - \mbf{S}\left(\mbg{\omega}(t)\right) \mbf{v}(t) + \mbf{g}(t)\\
			\dot{\mbf{g}}(t) = - \mbf{S}\left(\mbg{\omega}(t)\right) \mbf{g}(t)\\
			r_1(t) = \left\| \mbf{s}_1 - \mbf{p}(t) \right\|\\
			\vdots\\
			r_{n\low{L}}(t) = \left\| \mbf{s}_{n\low{L}} - \mbf{p}(t) \right\|
		\end{array}
	\right..
	\label{eq:PS:SD:SD}
\end{equation}
The problem considered in the paper is the design of a filter for
\eqref{eq:PS:SD:SD} assuming noisy measurements.

In the accelerometer measurements \eqref{eq:PS:SD:a} it was not
considered accelerometer bias. Its inclusion in the system dynamics
is trivial but the observability analysis would  result even more
tedious. Nevertheless, the observability of linear motion quantities
considering accelerometer bias was previously studied by the authors
in \cite{paper:Batista:ACC:2009:2} and it is rather straightforward
to conclude about the necessary (and sufficient) conditions for
observability of the system considering also accelerometer bias. For
the sake of simplicity and clarity of presentation, which focuses on
the novel Long Baseline solution, the accelerometer is assumed to be
calibrated in this paper and the overall system analysis and design
will be considered elsewhere.

\subsection{Long Baseline Configuration}
\label{sec:PS:LBLC}

Long Baseline acoustic configurations have been widely used in the
past in the design of navigation systems. In the remainder of the
paper the following assumption is considered:

\begin{assumption}
	There exists at least 4 noncoplanar landmarks.
	\label{assumption:PS:LBL Configuration}
\end{assumption}

When there exist at least 4 noncoplanar landmarks, it is always
possible to determine the position of the vehicle from the range
measurements. When there are fewer measurements that may not always
be possible. For example, for a static vehicle and 3 landmarks, 
there exist two possible solutions. Although it is assumed that
there are at least 4 noncoplanar landmarks, the proposed solution
is general and a more detailed discussion on different LBL
configurations will be presented in the sequel.

\section{Filter Design}
\label{sec:FD}

This section presents the main results of the paper. Firstly, a
state transformation is applied, in Section \ref{sec:FD:ST}, to
reduce the complexity of the system dynamics. Afterwards, state
augmentation is proposed, in Section \ref{sec:FD:SA}, in order to
derive sensor-based deterministic system dynamics, and in Section
\ref{sec:FD:OA} the observability of the resulting system is
analyzed. The design of a Kalman filter is detailed in Section
\ref{sec:FD:KF} and finally, in Section \ref{sec:FD:PC}, some
practical considerations are presented.

\subsection{State Transformation}
\label{sec:FD:ST}

Let
$\mbf{T}(t) := \rm{diag} \left( \mbf{I}, \mbf{R}(t), \mbf{R}(t) \right)$
and consider the state transformation
\begin{equation}
	\left[
		\begin{array}{c}
			\mbf{x}_1(t)\\
			\mbf{x}_2(t)\\
			\mbf{x}_3(t)
		\end{array}
	\right]
	:=
	\mbf{T}(t)
	\left[
		\begin{array}{c}
			\mbf{p}(t)\\
			\mbf{v}(t)\\
			\mbf{g}(t)
		\end{array}
	\right],
	\label{eq:FD:ST:ST}
\end{equation}
which is a Lyapunov state transformation previously used by the
authors, see \cite{paper:Batista:ACC:2009:1}. The new system
dynamics are given by
\begin{equation}
	\left\{
		\begin{array}{l}
			\dot{\mbf{x}}_1(t) = \mbf{x}_2(t)\\
			\dot{\mbf{x}}_2(t) = \mbf{x}_3(t)  + \mbf{u}(t)\\
			\dot{\mbf{x}}_3(t) = \mbf{0}\\
			r_1(t) = \left\| \mbf{s}_1 - \mbf{x}_1(t) \right\|\\
			\vdots\\
			r_{n\low{L}}(t) = \left\| \mbf{s}_{n\low{L}} - \mbf{x}_1(t) \right\|
		\end{array}
	\right.,
	\label{eq:FD:ST:SD}
\end{equation}
where $\mbf{u}(t) := \mbf{R}(t) \mbf{a}(t)$.
Notice that, as \eqref{eq:FD:ST:ST} is a Lyapunov state
transformation, all observability properties are preserved
\cite{book:Brockett:FiniteDimensionalLinearSystems}. The advantage
of considering the state transformation \eqref{eq:FD:ST:ST} is that
the new system is time invariant, although it is still nonlinear.

\subsection{State Augmentation}
\label{sec:FD:SA}

To derive a linear system that mimics the dynamics of the nonlinear
system \eqref{eq:FD:ST:SD}, define $n\low{L}+4$ additional scalar
state variables as
\[
	\left\{
		\begin{array}{l}
			x_4(t) := r_1(t)\\
			\vdots\\
			x_{n\low{L}+3}(t) := r_{n\low{L}}(t)\\
			x_{n\low{L}+4}(t) := \mbf{x}_1(t) \cdot \mbf{x}_2(t)\\
			x_{n\low{L}+5}(t) := \mbf{x}_1(t) \cdot \mbf{x}_3(t) + \left\| \mbf{x}_2(t) \right\|^2\\
			x_{n\low{L}+6} := \mbf{x}_2(t) \cdot \mbf{x}_3(t)\\
			x_{n\low{L}+7}(t) := \left\| \mbf{x}_3(t) \right\|^2
		\end{array}
	\right.
\]
and denote by
\[
	\mbg{x}(t) = \left[ \mbf{x}_1^T(t) \, \mbf{x}_2^T(t) \, \mbf{x}_3^T(t) \, x_4(t) \, \ldots \, x_{n\low{L}+7} \right]^T \in \mbb{R}^n,
\]
$n = 13 + n\low{L}$, the augmented state. It is easy to verify that
the augmented dynamics can be written as
\[
	\dot{\mbg{x}}(t) = \mbf{A}(t) \mbg{x}(t) + \mbf{B} \mbf{u}(t)\\
\]
where
\begin{eqnarray}
	& \mbf{A}(t) 
	=
	\left[\!
		\begin{array}{cccccccccc}
			\mbf{0} \! & \! \mbf{I} \! & \! \mbf{0} \! & \! \mbf{0} \! & \! \ldots \! & \! \mbf{0} \! & \! \mbf{0} \! & \! \mbf{0} \! & \! \mbf{0} \! & \! \mbf{0}\\
			\mbf{0} \! & \! \mbf{0} \! & \! \mbf{I} \! & \! \mbf{0} \! & \! \ldots \! & \! \mbf{0} \! & \! \mbf{0} \! & \! \mbf{0} \! & \! \mbf{0} \! & \! \mbf{0}\\
			\mbf{0} \! & \! \mbf{0} \! & \! \mbf{0} \! & \! \mbf{0} \! & \! \ldots \! & \! \mbf{0} \! & \! \mbf{0} \! & \! \mbf{0} \! & \! \mbf{0} \! & \! \mbf{0}\\
			\mbf{0} \! & \! -\frac{\mbf{s}_1^T}{r_1\left(t\right)} \! & \! \mbf{0} \! & \! 0 \! & \! \ldots \! & \! 0 \! & \!\frac{1}{r_1\left(t\right)} \! & \! 0 \! & \! 0 \! & \! 0\\
			\vdots \! & \! \vdots \! & \! \vdots \! & \! \vdots \! & \! \! & \! \vdots \! & \! \vdots \! & \! \vdots \! & \! \vdots \! & \! \vdots\\
			\mbf{0} \! & \! -\frac{\mbf{s}_{n_L}^T}{r_{n\low{L}}\left(t\right)} \! & \! \mbf{0} \! & \! 0 \! & \! \ldots \! & \! 0 \! & \! \frac{1}{r_{n_L}\left(t\right)} \! & \! 0 \! & \! 0 \! & \! 0\\
			\mbf{u}^T(t) \! & \! \mbf{0} \! & \! \mbf{0} \! & \! 0 \! & \! \ldots \! & \! 0 \! & \! 0 \! & \! 1 \! & \! 0 \! & \! 0\\
			\mbf{0} \! & \! 2 \mbf{u}^T(t) \! & \! \mbf{0} \! & \! 0 \! & \! \ldots \! & \! 0 \! & \! 0 \! & \! 0 \! & \! 3 \! & \! 0\\
			\mbf{0} \! & \!\mbf{0} \! & \! \mbf{u}^T(t) \! & \! 0 \! & \! \ldots \! & \! 0 \! & \! 0 \! & \! 0 \! & \! 0 \! & \! 1\\
			\mbf{0} \! & \! \mbf{0} \! & \! \mbf{0} \! & \! 0 \! & \! \ldots \! & \! 0 \! & \! 0 \! & \! 0 \! & \! 0 \! & \! 0
		\end{array}\!
	\right]
	&
	\label{eq:FD:SA:A}
\end{eqnarray}
and
$\mbf{B} = \left[ \mbf{0} \, \mbf{I} \, \mbf{0} \, 0 \, \ldots 0\right]^T$.

In order to complete the augmented system dynamics, notice first
that the states $x_4, \, \ldots, \, x_{3+n\low{L}}$ are measured. In
addition to that, it is straightforward to show that
\[
	\frac{2\left(\mbf{s}_i - \mbf{s}_j \right) \cdot \mbf{p}(t)}{r_i(t) + r_j(t)}  + r_i(t) - r_j(t) = \frac{\left\| \mbf{s}_i \right\|^2 - \left\| \mbf{s}_j \right\|^2}{r_i(t) + r_j(t)}
\]
or, equivalently,
\begin{equation}
	\frac{2\left(\mbf{s}_i - \mbf{s}_j \right) \cdot \mbf{x}_1(t)}{r_i(t) + r_j(t)}  + x_{3+i}(t) - x_{3+j}(t) = \frac{\left\| \mbf{s}_i \right\|^2 - \left\| \mbf{s}_j \right\|^2}{r_i(t) + r_j(t)},
	\label{eq:FD:SA:Additional Output}
\end{equation}
where the right side of \eqref{eq:FD:SA:Additional Output} is known
and the left side depends on the system state. Discarding the
original nonlinear system output, it is possible to write an
augmented system output as
{\small
\[
	\left\{ \!\!\!
		\begin{array}{l}
			y_1(t) = x_4(t)\\
			\vdots\\
			y_{n\low{L}}(t) = x_{3 + n\low{L}}(t)\\
			y_{n\low{L} + 1} = \frac{2\left(\mbf{s}_1 - \mbf{s}_2 \right) \cdot \mbf{x}_1(t)}{r_1(t) + r_2(t)}  + x_{3+1}(t) - x_{3+2}(t)\\
			y_{n\low{L} + 2} = \frac{2\left(\mbf{s}_1 - \mbf{s}_3 \right) \cdot \mbf{x}_1(t)}{r_1(t) + r_3(t)}  + x_{3+1}(t) - x_{3+3}(t)\\
			\vdots\\
			y_{n\low{L} + C_2^{n\low{L}}-1} = \frac{2\left(\mbf{s}_{n\low{L}-2} \! - \! \mbf{s}_{n\low{L}} \right) \! \cdot \! \mbf{x}_1(t)}{r_{n\low{L}-2}(t) + r_{n\low{L}}(t)}  + x_{3+{n\low{L}-2}}(t) \! - \! x_{3+{n\low{L}}}(t)\\
			y_{n\low{L} + C_2^{n\low{L}}} = \frac{2\left(\mbf{s}_{n\low{L}-1} \! - \! \mbf{s}_{n\low{L}} \right) \! \cdot \! \mbf{x}_1(t)}{r_{n\low{L}-1}(t) + r_{n\low{L}}(t)}  + x_{3+{n\low{L}-1}}(t) \! - \! x_{3+{n\low{L}}}(t)
		\end{array}
	\right.\!\!\!\!\!\!,
\]}\noindent
where $C_2^{n\low{L}}$ is the number of 2-combinations of $n\low{L}$
elements, i.e., $C_2^{n\low{L}} = \frac{\left( n\low{L} - 1\right) n\low{L}}{2}$.
In compact form, the augmented system dynamics are given by
\begin{equation}
	\left\{
		\begin{array}{l}
			\dot{\mbg{x}}(t) = \mbf{A}(t) \mbg{x}(t) + \mbf{B} \mbf{u}(t)\\
			\mbf{y}(t) = \mbf{C}(t) \mbg{x}(t)
		\end{array}
	\right.,
	\label{eq:FD:SA:SD}
\end{equation}
where 
\begin{equation}
	\mbf{C}(t)
	=
	\left[
		\begin{array}{ccc|c|c}
			\mbf{0}_{n\low{L}\times 3} & \mbf{0}_{n\low{L}\times 3} & \mbf{0}_{n\low{L}\times 3} & \mbf{I}_{n\low{L}} & \mbf{0}_{n\low{L}\times 4}\\
			\hline
			\mbf{C}_1(t) & \mbf{0}_{C_2^{n\low{L}}\times 3} & \mbf{0}_{C_2^{n\low{L}}\times 3} & \mbf{C}_2 & \mbf{0}_{C_2^{n\low{L}}\times 4}
		\end{array}
	\right],
	\label{eq:FD:SA:C}
\end{equation}
\[
	\mbf{C}_1(t) =
	\left[
		\begin{array}{c}
			\frac{2\left(\mbf{s}_1 - \mbf{s}_2 \right)^T}{r_1(t) + r_2(t)}\\
			\frac{2\left(\mbf{s}_1 - \mbf{s}_3 \right)^T}{r_1(t) + r_3(t)}\\
			\vdots\\
			\frac{2\left(\mbf{s}_{n\low{L}-2} - \mbf{s}_{n\low{L}} \right)^T}{r_{n\low{L}-2}(t) + r_{n\low{L}}(t)}\\
			\frac{2\left(\mbf{s}_{n\low{L}-1} - \mbf{s}_{n\low{L}} \right)^T}{r_{n\low{L}-1}(t) + r_{n\low{L}}(t)}
		\end{array}
	\right] \in \mbb{R}^{n\low{L} \times 3},
\]
and
\[
	\mbf{C_2}
	=
	\left[
		\begin{array}{ccccccc}
			1 & -1 & 0 & \ldots & \ldots & \ldots & 0\\
			1 & 0 & -1 & 0 & \ldots & \ldots & 0 \\
			& & & \vdots \\
			0 & \ldots & \ldots & 0 & 1 & 0 & -1\\
			0 & \ldots & \ldots & \ldots & 0 & 1 & -1
		\end{array}
	\right] \in \mbb{R}^{C_2^{n\low{L}}\times n\low{L}}
\]
is the matrix that encodes all the possible combinations of
differences of pairs of ranges.

The dynamic system \eqref{eq:FD:SA:SD} can be regarded as a linear
time-varying system, even though the system matrices $\mbf{A}(t)$
and $\mbf{C}(t)$ depend explicitly on the system input and output,
as evidenced by \eqref{eq:FD:SA:A} and \eqref{eq:FD:SA:C}.
However, these are known signals and therefore pose no problem,
other than the fact that the observability of the system may depend
on the system input, which does not happen for linear systems whose
system matrices do not depend on the system input. Also, notice that
there is nothing in \eqref{eq:FD:SA:SD} imposing
\begin{equation}
	\left\{
		\begin{array}{l}
			r_1(t) = \left\| \mbf{s}_1 - \mbf{x}_1(t) \right\|\\
			\vdots\\
			r_{n\low{L}}(t) = \left\| \mbf{s}_{n\low{L}} - \mbf{x}_1(t) \right\|\\
			x_{n\low{L}+4}(t) = \mbf{x}_1(t) \cdot \mbf{x}_2(t)\\
			x_{n\low{L}+5}(t) = \mbf{x}_1(t) \cdot \mbf{x}_3(t) + \left\| \mbf{x}_2(t) \right\|^2\\
			x_{n\low{L}+6}(t) = \mbf{x}_2(t) \cdot \mbf{x}_3(t)\\
			x_{n\low{L}+7}(t) = \left\| \mbf{x}_3(t) \right\|^2
		\end{array}
	\right..
	\label{eq:OA:AM:SA:Algebraic Restrictions}
\end{equation}
These restrictions could be easily implemented but this form is
preferred since it allows to consider the system as linear.
However, care must be taken when extrapolating conclusions from the
observability of \eqref{eq:FD:SA:SD} to the observability of
\eqref{eq:FD:ST:SD} or \eqref{eq:PS:SD:SD}. Finally, the following
assumption is required so that \eqref{eq:FD:SA:A} is well defined,
which is not restrictive from the practical point of view since it
would make no sense to have the vehicle at the same position of a
landmark, where an acoustic transponder is installed.
\begin{assumption}
	The motion of the vehicle is such that
	\[
		\begin{array}{cccc}
			\exists & \forall & : & r_m \leq r_i(t).\\
			r_m > 0 & t>t_0
		\end{array}
	\]
	\label{FD:SA:Assumption Ranges}
\end{assumption}

\subsection{Observability Analysis}
\label{sec:FD:OA}

In the previous section a LTV system was derived that aims to
capture the behavior of the original nonlinear system. The analysis
of the observability of this LTV system is carried out in this
section and its behavior compared with that of the nonlinear system.

In order to proceed with the analysis of the observability of the
LTV system \eqref{eq:FD:SA:SD}, it is convenient to compute the 
observability Gramian associated with the pair
$\left( \mbf{A}(t), \mbf{C}(t)\right)$ and, in order to do so, the
transition matrix associated with the system matrix $\mbf{A}(t)$.
Let
\[
	\mbf{u}^{[1]}\!\left(t,t\low{0}\right) \! := \int_{t\low{0}}^{t} \mbf{u}\left(\sigma\right) d\sigma
\]
and
\[
	\mbf{u}^{[2]}\!\left(t,t_0\right) \! := \int_{t\low{0}}^{t} \int_{t\low{0}}^{\sigma_1} \mbf{u}\left(\sigma_2\right) d\sigma_2 d\sigma_1.
\]
Long, but straightforward, computations show that the transition
matrix associated with $\mbf{A}(t)$ is given by
\[
	\mbg{\phi} \left(t, t\low{0}\right) =
	\left[
		\begin{array}{ccc}
			\mbg{\phi}\low{AA} \left(t, t\low{0}\right) & \mbf{0} & \mbf{0}\\
			\mbg{\phi}\low{BA} \left(t, t\low{0}\right) & \mbf{I} & \mbg{\phi}\low{BC} \left(t, t\low{0}\right)\\
			\mbg{\phi}\low{CA} \left(t, t\low{0}\right)  & \mbf{0} & \mbg{\phi}\low{CC} \left(t, t\low{0}\right)
		\end{array}
	\right],
\]
where
\[
	\mbg{\phi}\low{AA} \left(t, t\low{0}\right) =
	\left[
		\begin{array}{ccc}
			\mbf{I} & \left(t-t\low{0}\right) \mbf{I} & \frac{\left(t-t\low{0}\right)^2}{2} \mbf{I}\\
			\mbf{0} & \mbf{I} & \left(t-t\low{0}\right) \mbf{I}\\
			\mbf{0} & \mbf{0} & \mbf{I}
		\end{array}
	\right],
\]
\[
	\mbg{\phi}\low{BA} \left(t, t\low{0}\right)
	=
	\left[
		\begin{array}{ccc}
			\mbg{\phi}\low{BA1} \left(t, t\low{0}\right) & \mbg{\phi}\low{BA2} \left(t, t\low{0}\right) & \mbg{\phi}\low{BA3} \left(t, t\low{0}\right)
		\end{array}
	\right],
\]
\[
	\mbg{\phi}\low{BA1} \left(t, t\low{0}\right) =
	\left[ \!\!\!
		\begin{array}{c}
			\int_{t_0}^t \frac{\left[\mbf{u}^{[1]}\!\left(\sigma,t\low{0}\right)\right]^T}{r_1\left(\sigma\right)} d\sigma\\
			\vdots\\
			\int_{t_0}^t \frac{\left[\mbf{u}^{[1]}\!\left(\sigma,t\low{0}\right)\right]}{r_{n\low{L}}\left(\sigma\right)} d\sigma \\
		\end{array} \!\!\!
	\right]\!\!,
\]
\begin{eqnarray*}
	& \mbg{\phi}\low{BA2} \left(t, t\low{0}\right) = &
	\nonumber\\
	& \left[ \!\!\!
		\begin{array}{c}
			\int_{t_0}^t \frac{-\mbf{s}_1^T + \left(\sigma-t_0\right) \left[\mbf{u}^{[1]} \!\left(\sigma,t_0\right)\right]^T + \left[\mbf{u}^{[2]} \!\left(\sigma,t_0\right)\right]^T}{r_1\left(\sigma\right)} d\sigma
			\nonumber\\
			\vdots\\
			\int_{t_0}^t \frac{-\mbf{s}_{n\low{L}}^T + \left(\sigma-t_0\right) \left[\mbf{u}^{[1]} \!\left(\sigma,t_0\right)\right]^T + \left[\mbf{u}^{[2]} \!\left(\sigma,t_0\right)\right]^T}{r_{n\low{L}}\left(\sigma\right)} d\sigma
		\end{array} \!\!\!
	\right]\!\!, &
\end{eqnarray*}
{\small
\begin{eqnarray*}
	& \mbg{\phi}\low{BA3} \left(t, t\low{0}\right) = &
	\nonumber\\
	&
	\left[ \!\!\!
		\begin{array}{c}
			\int_{t_0}^t \!\! \left(\sigma-t\low{0}\right) \frac{ - \mbf{s}_1^T  + \frac{\left(\sigma-t_0\right)}{2} \left[\mbf{u}^{[1]} \!\left(\sigma,t_0\right)\right]^T + \left[\mbf{u}^{[2]} \!\left(\sigma,t_0\right)\right]^T}{r_1\left(\sigma\right)} d\sigma \\
			\vdots\\
			\int_{t_0}^t \!\! \left(\sigma-t\low{0}\right) \frac{ - \mbf{s}_{n\low{L}}^T  + \frac{\left(\sigma-t_0\right)}{2} \left[\mbf{u}^{[1]} \!\left(\sigma,t_0\right)\right]^T + \left[\mbf{u}^{[2]} \!\left(\sigma,t_0\right)\right]^T}{r_{n\low{L}}\left(\sigma\right)} d\sigma \\
		\end{array} \!\!\!
	\right]\!\!,
	&
\end{eqnarray*}
\begin{eqnarray*}
	& \mbg{\phi}\low{BC} \left(t, t\low{0}\right) = &
	\nonumber\\
	&
	\!\!\!\left[ \!\!\!
		\begin{array}{cccc}
			\int_{t_0}^t \! \frac{1}{r_1\left(\sigma\right)} d\sigma \! & \! \int_{t_0}^t \! \frac{\sigma-t\low{0}}{r_1\left(\sigma\right)} d\sigma \! & \! \int_{t_0}^t \! \frac{3}{2} \frac{\left(\sigma-t\low{0}\right)^2}{r_1\left(\sigma\right)} d\sigma \! & \! \int_{t_0}^t \! \frac{1}{2} \frac{\left(\sigma-t\low{0}\right)^3}{r_1\left(\sigma\right)} d\sigma\\
			\vdots & \vdots & \vdots & \vdots\\
			\int_{t_0}^t \! \frac{1}{r_{n\low{L}}\!\left(\sigma\right)} d\sigma \!\! & \!\! \int_{t_0}^t \! \frac{\sigma-t\low{0}}{r_{n\low{L}}\!\left(\sigma\right)} d\sigma \!\! & \!\! \int_{t_0}^t \! \frac{3}{2} \frac{\left(\sigma-t\low{0}\right)^2}{r_{n\low{L}}\!\left(\sigma\right)} d\sigma \!\! & \!\! \int_{t_0}^t \! \frac{1}{2} \frac{\left(\sigma-t\low{0}\right)^3}{r_{n\low{L}}\!\left(\sigma\right)} d\sigma\\
		\end{array} \!\!\!
	\right]\!\!,
	&
\end{eqnarray*}}\noindent
and $\mbg{\phi}\low{CA} \left(t, t\low{0}\right)$ and
$\mbg{\phi}\low{CC} \left(t, t\low{0}\right)$ are omitted for the
sake of simplicity, as they are not required in the sequel. The
observability Gramian associated with the pair
$\left( \mbf{A}(t), \mbf{C}(t)\right)$ is simply given by
\[
	\mbg{\cal W}\left(t_0, t_f\right) = \int_{t_0}^{t_f} {\mbg{\phi}^T\left(t, t_0\right) \mbf{C}^T(t) \mbf{C}(t) \mbg{\phi}\left(t, t_0\right) dt}.
\]

The following theorem establishes the observability of the LTV
system \eqref{eq:FD:SA:SD}.

\begin{theorem}
\label{thm:FD:OA:Obs LTV}
The linear time-varying system \eqref{eq:FD:SA:SD} is observable on
$\left[t_0, t_f\right]$, $t_0 < t_f$.
\end{theorem}
\begin{pf}
The proof follows by contradiction. Suppose that the LTV system
\eqref{eq:FD:SA:SD} is not observable on
${\cal I} := \left[t_0, t_f\right]$. Then, there exists a nonnull
vector $\mbf{d} \in \mbb{R}^{13+n\low{L}}$
\[
	\mbf{d} =
	\left[
			\mbf{d}_1^T
			\mbf{d}_2^T
			\mbf{d}_3^T
			\mbf{d}_4^T
			d_5
			d_6
			d_7
			d_8
	\right]^T,
\]
with $\mbf{d}_1, \, \mbf{d}_2, \,\mbf{d}_3 \in \mbb{R}^3$,
$\mbf{d}_4 \in \mbb{R}^{n\low{L}}$,
$ d_5, \,\ldots, \,d_8 \in \mbb{R}$, such that
$\mbf{d}^T \mbg{\cal W}\left(t_0, t\right) \mbf{d} = 0$
for all $t \in {\cal I}$ or, equivalently,
\begin{equation}
	\int_{t\low{0}}^t \left\| \mbf{C}\left(\tau\right) \mbg{\phi} \left(\tau, t\low{0}\right) \mbf{d} \right\|^2 d\tau = 0.
	\label{eq:FD:OA:OBS L:S1}
\end{equation}
Taking the time derivative of \eqref{eq:FD:OA:OBS L:S1} yields
\[
	\left\| \mbf{C}\left(t\right) \mbg{\phi} \left(t, t\low{0}\right) \mbf{d} \right\|^2 = 0,
\]
which implies, in particular, that
\begin{equation}
	\mbf{C}\left(t\right) \mbg{\phi} \left(t, t\low{0}\right) \mbf{d} = 0
	\label{eq:FD:OA:OBS L:S2}
\end{equation}
for all $t \in {\cal I}$. Substituting $t = t\low{0}$ in
\eqref{eq:FD:OA:OBS L:S2} gives
\begin{equation}
	\left[
		\begin{array}{c}
			\mbf{d}_4\\
			\mbf{C}_1\left(t\low{0}\right) \mbf{d}_1 + \mbf{C}_2 \mbf{d}_4
		\end{array}
	\right]
	= \mbf{0}
	\label{eq:FD:OA:OBS L:S3}
\end{equation}
which implies immediately that
\begin{equation}
	\mbf{d}_4 = \mbf{0}.
	\label{eq:FD:OA:OBS L:S4}
\end{equation}
Substituting \eqref{eq:FD:OA:OBS L:S4} in \eqref{eq:FD:OA:OBS L:S3}
gives
\begin{equation}
	\left[
		\begin{array}{c}
			\frac{2}{r_1\left(t_0\right) + r_2\left(t_0\right)} \left( \mbf{s}_1 - \mbf{s}_2\right)^T\\
			\frac{2}{r_1\left(t_0\right) + r_3\left(t_0\right)} \left( \mbf{s}_1 - \mbf{s}_3\right)^T\\
			\vdots\\
			\frac{2}{r_{n\low{L}-2}\left(t_0\right) + r_{n\low{L}}\left(t_0\right)} \left( \mbf{s}_{n\low{L}-2} - \mbf{s}_{n\low{L}}\right)^T\\
			\frac{2}{r_{n\low{L}-1}\left(t_0\right) + r_{n\low{L}}\left(t_0\right)} \left( \mbf{s}_{n\low{L}-1} - \mbf{s}_{n\low{L}}\right)^T
		\end{array}
	\right]
	\mbf{d}_1
	= \mbf{0}.
	\label{eq:FD:OA:OBS L:S5}
\end{equation}
It is straightforward to show that, under Assumption
\ref{assumption:PS:LBL Configuration}, the only solution of
\eqref{eq:FD:OA:OBS L:S5} is
\begin{equation}
	\mbf{d}_1 = \mbf{0}.
	\label{eq:FD:OA:OBS L:S6}
\end{equation}
Now, from \eqref{eq:FD:OA:OBS L:S2}, it is possible to write
\begin{equation}
	\frac{d}{dt} \left[\mbf{C}\left(t\right) \mbg{\phi} \left(t, t\low{0}\right) \mbf{d} \right]= \mbf{0}
	\label{eq:FD:OA:OBS L:S7}
\end{equation}
for all $t \in {\cal I}$. Expanding \eqref{eq:FD:OA:OBS L:S7}, and
considering \eqref{eq:FD:OA:OBS L:S4} and \eqref{eq:FD:OA:OBS L:S6},
allows to write
{\small
\begin{eqnarray*}
	& \frac{-\mbf{s}_i^T \mbf{d}_2 + \left(t-t_0\right) \left[\mbf{u}^{[1]} \!\left(t,t_0\right)\right]^T \mbf{d}_2 + \left[\mbf{u}^{[2]} \!\left(t,t_0\right)\right]^T \mbf{d}_2}{r_i\left(t\right)} &
	\nonumber\\
	& + \left(t-t\low{0}\right)\frac{ - \mbf{s}_i^T \mbf{d}_3  + \frac{\left(t-t_0\right)}{2} \left[\mbf{u}^{[1]} \!\left(t,t_0\right)\right]^T \mbf{d}_3 +  \left[\mbf{u}^{[2]} \!\left(t,t_0\right)\right]^T \mbf{d}_3}{r_i\left(t\right)} &
	\nonumber\\
	& + \frac{1}{r_i\left(t\right)} d_5 + \frac{t-t\low{0}}{r_i\left(t\right)} d_6 + \frac{3}{2} \frac{\left(t-t\low{0}\right)^2}{r_i\left(t\right)} d_7 + \frac{1}{2} \frac{\left(t-t\low{0}\right)^3}{r_i\left(t\right)} d_8 = 0 &
\end{eqnarray*}}\noindent
for all $i = 1, \, \ldots, \, n\low{L}$, $t \in I$ or,
equivalently,
\begin{eqnarray}
	& \left. -\mbf{s}_i^T \mbf{d}_2 + \left(t-t_0\right) \left[\mbf{u}^{[1]} \!\left(t,t_0\right)\right]^T \mbf{d}_2 + \left[\mbf{u}^{[2]} \!\left(t,t_0\right)\right]^T \mbf{d}_2\right. &
	\nonumber\\
	& \left. -\left(t-t\low{0}\right) \mbf{s}_i^T \mbf{d}_3  + \frac{\left(t-t_0\right)^2}{2} \left[\mbf{u}^{[1]} \!\left(t,t_0\right)\right]^T \mbf{d}_3 \right. &
	\nonumber\\
	& + \left(t-t_0\right) \left[\mbf{u}^{[2]} \!\left(t,t_0\right)\right]^T \mbf{d}_3 + d_5 + \left(t-t\low{0}\right) d_6 &
	\nonumber\\
	&  + \frac{3}{2} \left(t-t\low{0}\right)^2 d_7 + \frac{1}{2} \left(t-t\low{0}\right)^3 d_8 = 0. &
	\label{eq:FD:OA:OBS L:S8}
\end{eqnarray}
Substituting $t = t_0$ in \eqref{eq:FD:OA:OBS L:S8} gives
\begin{equation}
	\left[
		\begin{array}{cc}
			-\mbf{s}_1^T & 1 \\
			-\mbf{s}_2^T & 1 \\
			\vdots\\
			-\mbf{s}_{n\low{L}-1}^T & 1 \\
			-\mbf{s}_{n\low{L}}^T & 1
		\end{array}
	\right]
	\left[
		\begin{array}{c}
			\mbf{d}_2\\
			d_5
		\end{array}
	\right]
	= \mbf{0}.
	\label{eq:FD:OA:OBS L:S9}
\end{equation}
Again, it is straightforward to show that, under Assumption
\ref{assumption:PS:LBL Configuration}, the only solution of
\eqref{eq:FD:OA:OBS L:S9} is
\begin{equation}
	\left\{
		\begin{array}{l}
			\mbf{d}_2 = \mbf{0}\\
			d_5 = 0
		\end{array}
	\right..
	\label{eq:FD:OA:OBS L:S10}
\end{equation}
Now, considering \eqref{eq:FD:OA:OBS L:S10} in
\eqref{eq:FD:OA:OBS L:S8} and taking its time derivative gives
\begin{eqnarray}
	& \left. - \mbf{s}_i^T \mbf{d}_3  + \frac{\left(t-t_0\right)^2}{2} \left[\mbf{u}\left(t\right)\right]^T \mbf{d}_3 + \left(t-t_0\right) \left[\mbf{u}^{[1]} \!\left(t,t_0\right)\right]^T \mbf{d}_3 \right. &
	\nonumber\\
	& + \left(t-t_0\right) \left[\mbf{u}^{[1]} \!\left(t,t_0\right)\right]^T \mbf{d}_3 + \left[\mbf{u}^{[2]} \!\left(t,t_0\right)\right]^T \mbf{d}_3 + d_6 &
	\nonumber\\
	&  + 3 \left(t-t\low{0}\right) d_7 + \frac{3}{2} \left(t-t\low{0}\right)^2 d_8 = 0 &
	\label{eq:FD:OA:OBS L:S11}
\end{eqnarray}
for all $i = 1, \, \ldots, \, n\low{L}$ and $t \in I$.
With $t = t_0$ in \eqref{eq:FD:OA:OBS L:S11},
\begin{equation}
	\left[
		\begin{array}{cc}
			-\mbf{s}_1^T & 1 \\
			-\mbf{s}_2^T & 1 \\
			\vdots\\
			-\mbf{s}_{n\low{L}-1}^T & 1 \\
			-\mbf{s}_{n\low{L}}^T & 1
		\end{array}
	\right]
	\left[
		\begin{array}{c}
			\mbf{d}_3\\
			d_6
		\end{array}
	\right]
	= \mbf{0}.
	\label{eq:FD:OA:OBS L:S12}
\end{equation}
Again, under Assumption \ref{assumption:PS:LBL Configuration}, the
only solution of \eqref{eq:FD:OA:OBS L:S12} is
\begin{equation}
	\left\{
		\begin{array}{l}
			\mbf{d}_3 = \mbf{0}\\
			d_6 = 0
		\end{array}
	\right..
	\label{eq:FD:OA:OBS L:S13}
\end{equation}
Finally, substituting \eqref{eq:FD:OA:OBS L:S13} in
\eqref{eq:FD:OA:OBS L:S11} gives
\begin{eqnarray}
	& 3 \left(t-t\low{0}\right) d_7 + \frac{3}{2} \left(t-t\low{0}\right)^2 d_8 = 0 &
	\label{eq:FD:OA:OBS L:S14}
\end{eqnarray}
for all $t \in I$. Since $t-t_0$ and $\left(t-t_0\right)^2$ are
linearly independent functions, the only solution of
\eqref{eq:FD:OA:OBS L:S14} is $d_7 = d_8 = 0$. But that contradicts 
the hypothesis that there exists a nonnull vector $\mbf{d}$ such
that \eqref{eq:FD:OA:OBS L:S1} is true. Therefore, the LTV system
\eqref{eq:FD:SA:SD} is observable.
\end{pf}

The fact that \eqref{eq:FD:SA:SD} is observable does not mean that
the nonlinear system \eqref{eq:FD:ST:SD} is observable nor that an
observer for \eqref{eq:FD:SA:SD} is also an observer for
\eqref{eq:FD:ST:SD}, as there is nothing in the system dynamics 
\eqref{eq:FD:SA:SD} imposing the algebraic restrictions
\eqref{eq:OA:AM:SA:Algebraic Restrictions}. However, this turns out
to be true, as shown in the following theorem.

\begin{theorem}
\label{thm:FD:OA:Obs NS}
The nonlinear system \eqref{eq:FD:ST:SD} is observable in the sense
that, given the output
$\left\{ \mbf{y}(t), \, t \in \left[t_0, t_f \right] \right\}$ and
the input
$\left\{ \mbf{u}(t), \, t \in \left[t_0, t_f \right] \right\}$, the
initial state
$\mbf{x}\left(t_0\right) = \left[\mbf{x}_1^T\left(t_0\right) \, \mbf{x}_2^T\left(t_0\right) \, \mbf{x}_3^T\left(t_0\right) \right]^T$
is uniquely defined. Moreover, a state observer for the LTV system
\eqref{eq:FD:SA:SD} with globally asymptotically stable error
dynamics is also a state observer for the nonlinear system
\eqref{eq:FD:ST:SD}, with globally asymptotically stable error
dynamics.
\end{theorem}
\begin{pf}
It has been shown, in Theorem \ref{thm:FD:OA:Obs LTV}, that the LTV
system \eqref{eq:FD:SA:SD} is observable. Therefore, given the output 
$\left\{ \mbf{y}(t), \, t \in \left[t_0, t_f \right] \right\}$ and
the input
$\left\{ \mbf{u}(t), \, t \in \left[t_0, t_f \right] \right\}$, the
initial state of \eqref{eq:FD:SA:SD} is uniquely defined. Let
\[
	\mbf{z}\left(t_0\right) =
	\left[
			\mbf{z}_1^T\left(t_0\right)
			\mbf{z}_2^T\left(t_0\right)
			\mbf{z}_3^T\left(t_0\right)
			z_4\left(t_0\right)
			\ldots\\
			z_{n\low{L} + 7}\left(t_0\right)
	\right]^T
\]
$\mbf{z}_1, \, \mbf{z}_2, \,\mbf{z}_3 \in \mbb{R}^3$,
$z_4, \,\ldots, \,z_{n\low{L} + 7}\left(t_0\right) \in \mbb{R}$,
be the initial state of \eqref{eq:FD:SA:SD} and
$
	\mbf{x}\left(t_0\right) =
		\left[
				\mbf{x}_1^T\left(t_0\right)
				\mbf{x}_2^T\left(t_0\right)
				\mbf{x}_3^T\left(t_0\right)
		\right]^T
$
be the initial state of the nonlinear system \eqref{eq:FD:ST:SD}.
It is easy to show that the evolution, for the nonlinear system,
of $\mbf{x}_1(t)$, is given by
\begin{eqnarray}
	\mbf{x}_1(t) & = & \mbf{x}_1\left(t_0\right) + \left(t-t_0\right) \mbf{x}_2\left(t_0\right) + \frac{\left(t-t_0\right)^2}{2} \mbf{x}_3\left(t_0\right)
	\nonumber\\
	&& + \mbf{u}^{[2]} \!\left(t,t_0\right),
	\label{eq:FD:OA:OBS NL:x_1 NL}
\end{eqnarray}
while the output of the nonlinear system satisfies
\begin{eqnarray}
	& r_i^2(t) = \left\| \mbf{x}_1\left(t_0\right) - \mbf{s}_i \right\|^2 + \left(t-t_0\right)^2 \left\| \mbf{x}_2\left(t_0\right) \right\|^2 &
	\nonumber\\
	& + \frac{\left(t-t_0\right)^4}{4} \left\| \mbf{x}_3\left(t_0\right) \right\|^2 + 2 \left(t-t_0\right) \mbf{x}_1\left(t_0\right) \cdot \mbf{x}_2\left(t_0\right) & 
	\nonumber\\
	& + \left(t-t_0\right)^2 \mbf{x}_1\left(t_0\right) \cdot \mbf{x}_3\left(t_0\right) + \left(t-t_0\right)^3 \mbf{x}_2\left(t_0\right) \cdot \mbf{x}_3\left(t_0\right) &
	\nonumber\\
	& - 2 \left(t-t_0\right) \mbf{s}_i \cdot \mbf{x}_2\left(t_0\right) - \left(t-t_0\right)^2 \mbf{s}_i \cdot \mbf{x}_3\left(t_0\right) &
	\nonumber\\
	& + 2 \mbf{x}_1\left(t_0\right) \cdot \mbf{u}^{[2]} \!\left(t,t_0\right) + 2 \mbf{x}_2\left(t_0\right) \cdot \left(t-t_0\right) \mbf{u}^{[2]} \!\left(t,t_0\right) &
	\nonumber\\
	& + \mbf{x}_3\left(t_0\right) \cdot \left(t-t_0\right)^2 \mbf{u}^{[2]} \!\left(t,t_0\right) + \left\| \mbf{u}^{[2]} \!\left(t,t_0\right) \right\|^2 &
	\nonumber\\
	& - 2 \mbf{s}_i \cdot \mbf{u}^{[2]} \!\left(t,t_0\right). &
	\label{eq:FD:OA:OBS NL:r_i NL}
\end{eqnarray}
Therefore, from \eqref{eq:FD:OA:OBS NL:x_1 NL} and
\eqref{eq:FD:OA:OBS NL:r_i NL}, it is trivially shown that
\begin{eqnarray}
	& r_i^2(t) - r_j^2(t) + 2 \left( \mbf{s}_i - \mbf{s}_j \right) \cdot \mbf{x}_1(t) = \left\| \mbf{x}_1\left(t_0\right) - \mbf{s}_i \right\|^2 &
	\nonumber\\
	& - \left\| \mbf{x}_1\left(t_0\right) - \mbf{s}_j \right\|^2 + 2 \left( \mbf{s}_i - \mbf{s}_j \right) \cdot \mbf{x}_1\left(t_0\right) &
	\label{eq:OA:FD:OBS NL:S1}
\end{eqnarray}
for all $i,j \in \left\{1, \,\ldots, \,n\low{L} \right\}$ and
$t \in {\cal I}$.
Now, notice that, for the LTV system, multiplying the set of
augmented outputs 
$y_{n\low{L}+1}, \, \ldots,\, y_{n\low{L} + C_2^{n\low{L}}}$
by the corresponding sums of pair of ranges yields
\[
	r_i^2(t) - r_j^2(t) + 2 \left( \mbf{s}_i - \mbf{s}_j \right) \cdot \mbf{x}_1(t),
\]
while its evolution can be shown to satisfy
\begin{eqnarray}
	& r_i^2(t) - r_j^2(t) + 2 \left( \mbf{s}_i - \mbf{s}_j \right) \cdot \mbf{x}_1(t) = z_{3+i}^2\left(t_0\right) - z_{3+j}^2\left(t_0\right) &
	\nonumber\\
	& + 2 \left( \mbf{s}_i - \mbf{s}_j \right) \cdot \mbf{z}_1\left(t_0\right) &
	\label{eq:OA:FD:OBS NL:S2}
\end{eqnarray}
for all $i,j \in \left\{1, \,\ldots, \,n\low{L} \right\}$ and
$t \in {\cal I}$. The states of the augmented LTV system
$x_{4}(t), \, \dots, \, x_{n\low{L} + 7}(t)$ are actually measured
and correspond to the range measurements. Therefore, it must be
\begin{equation}
	z_{3+i}\left(t_0\right) = \left\| \mbf{x}_1\left(t_0\right) - \mbf{s}_i \right\|
	\label{eq:FD:OA:OBS NL:S3}
\end{equation}
for all $i \in \left\{1, \,\ldots, \,n\low{L} \right\}$. Then,
from the comparison of \eqref{eq:OA:FD:OBS NL:S1} and
\eqref{eq:OA:FD:OBS NL:S2}, and considering
\eqref{eq:FD:OA:OBS NL:S3}, it follows that
\begin{equation}
	\left[
		\begin{array}{c}
			\left( \mbf{s}_1 - \mbf{s}_2\right)^T\\
			\left( \mbf{s}_1 - \mbf{s}_3\right)^T\\
			\vdots\\
			\left( \mbf{s}_{n\low{L}-2} - \mbf{s}_{n\low{L}}\right)^T\\
			\left( \mbf{s}_{n\low{L}-1} - \mbf{s}_{n\low{L}}\right)^T
		\end{array}
	\right]
	\left[ \mbf{x}_1\left(t_0\right) - \mbf{z}_1\left(t_0\right) \right]
	= \mbf{0}.
	\label{eq:FD:OA:OBS NL:S4}
\end{equation}
Under Assumption \ref{assumption:PS:LBL Configuration} the only
solution of \eqref{eq:FD:OA:OBS NL:S4} is
\[
	\mbf{x}_1\left(t_0\right) = \mbf{z}_1\left(t_0\right).
\]
From \eqref{eq:FD:OA:OBS NL:r_i NL} it is possible to write
\begin{eqnarray}
	& r_i^2(t) - r_j^2(t) = \left\| \mbf{x}_1\left(t_0\right) - \mbf{s}_i \right\|^2 - \left\| \mbf{x}_1\left(t_0\right) - \mbf{s}_j \right\|^2 &
	\nonumber\\
	& - 2 \left(t-t_0\right) \left( \mbf{s}_i - \mbf{s}_j \right)\cdot \mbf{x}_2\left(t_0\right) &
	\nonumber\\
	& - \left(t-t_0\right)^2 \left( \mbf{s}_i - \mbf{s}_j \right) \cdot \mbf{x}_3\left(t_0\right) & 
	\nonumber\\
	& - 2 \left( \mbf{s}_i - \mbf{s}_j \right) \cdot \mbf{u}^{[2]} \!\left(t,t_0\right) &
	\label{eq:FD:OA:OBS NL:S5}
\end{eqnarray}
for all $i,j \in \left\{1, \,\ldots, \,n\low{L} \right\}$ and
$t \in {\cal I}$. On the other hand, the evolution of the square of
range for the LTV system \eqref{eq:FD:SA:SD} can be written as
\begin{eqnarray*}
	& r_i^2(t) = 2 \left[ \mbf{z}_1\left(t_0\right) - \mbf{s}_i \right] \cdot \mbf{u}^{[2]} \!\left(t,t_0\right) &
	\nonumber\\
	& - 2 \left(t-t_0\right) \mbf{s}_i \cdot \mbf{z}_2\left(t_0\right) + 2 \left(t-t_0\right) \mbf{u}^{[2]} \!\left(t,t_0\right) \cdot \mbf{z}_2\left(t_0\right) &
	\nonumber\\
	& - \left(t-t_0\right)^2 \mbf{s}_i \cdot \mbf{z}_3\left(t_0\right) + \left(t-t_0\right)^2 \mbf{u}^{[2]} \!\left(t,t_0\right) \cdot \mbf{z}_3\left(t_0\right) &
	\nonumber\\
	& + z_{3+i}^2\left(t_0\right) + 2 \left(t-t_0\right) z_{n\low{L} + 4}\left(t_0\right) &
	\nonumber\\
	& + \left(t-t_0\right)^2 z_{n\low{L} + 5}\left(t_0\right) + \left(t-t_0\right)^3 z_{n\low{L} + 6}\left(t_0\right) &
	\nonumber\\
	& + \frac{\left(t-t_0\right)^3}{4} z_{n\low{L} + 7}\left(t_0\right) + \left\| \mbf{u}^{[2]} \!\left(t,t_0\right) \right\|^2 &
\end{eqnarray*}
for all $i \in \left\{1, \,\ldots, \,n\low{L} \right\}$. Therefore,
it is possible to write
\begin{eqnarray}
	& r_i^2(t) - r_j^2(t) = z_{3+i}^2\left(t_0\right) - z_{3+j}^2\left(t_0\right) &
	\nonumber\\
	& - 2 \left(t-t_0\right) \left( \mbf{s}_i - \mbf{s}_j \right) \cdot \mbf{z}_2\left(t_0\right) &
	\nonumber\\
	& - \left(t-t_0\right)^2 \left( \mbf{s}_i - \mbf{s}_j \right) \cdot \mbf{z}_3\left(t_0\right) &
	\nonumber\\
	& - 2 \left( \mbf{s}_i - \mbf{s}_j \right) \cdot \mbf{u}^{[2]} \!\left(t,t_0\right) &
	\label{eq:FD:OA:OBS NL:S6}
\end{eqnarray}
for all $i,j \in \left\{1, \,\ldots, \,n\low{L} \right\}$ and
$t \in {\cal I}$. Taking the time derivative of both
\eqref{eq:FD:OA:OBS NL:S5} and \eqref{eq:FD:OA:OBS NL:S6} and
comparing for $t = t_0$ gives
\begin{equation}
	\left[
		\begin{array}{c}
			\left( \mbf{s}_1 - \mbf{s}_2\right)^T\\
			\left( \mbf{s}_1 - \mbf{s}_3\right)^T\\
			\vdots\\
			\left( \mbf{s}_{n\low{L}-2} - \mbf{s}_{n\low{L}}\right)^T\\
			\left( \mbf{s}_{n\low{L}-1} - \mbf{s}_{n\low{L}}\right)^T
		\end{array}
	\right]
	\left[ \mbf{x}_2\left(t_0\right) - \mbf{z}_2\left(t_0\right) \right]
	= \mbf{0}.
	\label{eq:FD:OA:OBS NL:S7}
\end{equation}
Again, under Assumption \ref{assumption:PS:LBL Configuration} the only
solution of \eqref{eq:FD:OA:OBS NL:S7} is
\[
	\mbf{x}_2\left(t_0\right) = \mbf{z}_2\left(t_0\right).
\]
Finally, taking the second time derivative of both 
\eqref{eq:FD:OA:OBS NL:S5} and \eqref{eq:FD:OA:OBS NL:S6} and
comparing for $t = t_0$ gives
\begin{equation}
	\left[
		\begin{array}{c}
			\left( \mbf{s}_1 - \mbf{s}_2\right)^T\\
			\left( \mbf{s}_1 - \mbf{s}_3\right)^T\\
			\vdots\\
			\left( \mbf{s}_{n\low{L}-2} - \mbf{s}_{n\low{L}}\right)^T\\
			\left( \mbf{s}_{n\low{L}-1} - \mbf{s}_{n\low{L}}\right)^T
		\end{array}
	\right]
	\left[ \mbf{x}_3\left(t_0\right) - \mbf{z}_3\left(t_0\right) \right]
	= \mbf{0}.
	\label{eq:FD:OA:OBS NL:S8}
\end{equation}
Again, under Assumption \ref{assumption:PS:LBL Configuration} the only
solution of \eqref{eq:FD:OA:OBS NL:S8} is
$\mbf{x}_3\left(t_0\right) = \mbf{z}_3\left(t_0\right)$.
This concludes the proof, as the initial state of the LTV system
\eqref{eq:FD:SA:SD}, which is uniquely defined, matches the initial
state of the nonlinear system \eqref{eq:FD:ST:SD}.
\end{pf}

\begin{remark}
	The concept of observability for nonlinear systems is not as strong
	as that presented in the statement of Theorem
	\ref{thm:FD:OA:Obs NS}, see
	\cite{paper:Krener:NonlinearControllabilityObservability}. That is
	the reasoning behind explicitly describing in what sense the system
	is observable.
\end{remark}

\begin{remark}
	In the proof of Theorem \ref{thm:FD:OA:Obs NS} it was not shown
	that all the algebraic relations
	\eqref{eq:OA:AM:SA:Algebraic Restrictions} are satisfied. It was
	only shown that the initial state of the LTV system
	\eqref{eq:FD:SA:SD} coincides with the initial state of the
	nonlinear system \eqref{eq:FD:ST:SD}. However, it is trivial to 
	show that those relations are indeed preserved.
\end{remark}

\begin{remark}
	Before concluding this section, it is important to remark that,
	although the observability results were derived with respect to the
	nonlinear system \eqref{eq:FD:ST:SD}, they also apply to the original
	nonlinear system \eqref{eq:PS:SD:SD} as they are related through a
	Lyapunov transformation. Also, the design of an observer for the
	original nonlinear system follows simply by reversing the state
	transformation \eqref{eq:FD:ST:ST}, as it will be detailed in the
	following section.
\end{remark}

\subsection{Kalman Filter}
\label{sec:FD:KF}

Although all the results derived so far were presented in a
deterministic setting, in practice there exists measurement noise
and often system disturbances. Therefore, a filtering solution is
proposed in this section instead of an observer. On the other hand,
Theorem \ref{thm:FD:OA:Obs NS} provides a constructive result in the
sense that a dynamic system with globally asymptotically stable error
dynamics for the LTV system \eqref{eq:FD:ST:SD} provides globally
asymptotically stable error dynamics for the estimation of the
state of the nonlinear system. Therefore, the design of a Kalman
filter follows for the LTV system \eqref{eq:FD:SA:SD}, albeit other
solutions could be devised, e.g., an ${\cal H}_\infty$ filter. It
is important to stress, however, that this filter is not optimal.
Indeed, looking into the system matrices, it is easy to see that, in
the presence of noise in the range or acceleration measurements,
there exists multiplicative noise. Nevertheless, the Kalman filter
is GAS, as it is straightforward to show that the system under
consideration is not only observable but also uniformly completely
observable.

In order to recover the augmented system dynamics in the original
coordinate space, consider the augmented state transformation
$\mbg{\chi}(t) = \mbf{T}_c(t) \mbg{x}(t)$,
where $\mbf{T}_c(t) = \rm{diag}\left(\mbf{I}, \mbf{R}^T(t), \mbf{R}^T(t), 1, \ldots, 1 \right)$.
Then, the nominal augmented system dynamics in the original
coordinate space are given by
\[
	\left\{
		\begin{array}{l}
			\dot{\mbg{\chi}}(t) = \mbg{\cal A}(t) \mbg{\chi}(t) + \mbf{B}(t) \mbf{a}(t)\\
			\mbf{y}(t) = \mbf{C}(t) \mbg{\chi}(t)
		\end{array}
	\right.,
\]
where
\begin{eqnarray*}
	& \mbg{\cal A}(t) = 
	\!\!\!\left[\!
		\begin{array}{cccccccc}
			\mbf{0} \! & \! \mbf{R}(t) \! & \! \mbf{0} \! & \! \mbf{0} \! & \! \mbf{0} \! & \! \mbf{0} \! & \! \mbf{0} \! & \! \mbf{0}\\
			\mbf{0} \! & \! - \mbf{S}\left(\mbg{\omega}(t)\right) \! & \! \mbf{I} \! & \! \mbf{0} \! & \! \mbf{0} \! & \! \mbf{0} \! & \! \mbf{0} \! & \! \mbf{0}\\
			\mbf{0} \! & \! \mbf{0} \! & \! - \mbf{S}\left(\mbg{\omega}(t)\right) \! & \! \mbf{0} \! & \! \mbf{0} \! & \! \mbf{0} \! & \! \mbf{0} \! & \! \mbf{0}\\
			\mbf{0} \! & \! -\frac{\mbf{s}_1^T \mbf{R}(t)}{r_1\left(t\right)} \! & \! \mbf{0} \! & \! 0 \! & \! \frac{1}{r_1\left(t\right)} \! & \! 0 \! & \! 0 \! & \! 0\\
			\vdots \! & \! \vdots \! & \! \vdots \! & \! \vdots \! & \! \vdots \! & \! \vdots \! & \! \vdots \! & \! \vdots\\
			\mbf{0} \! & \! -\frac{\mbf{s}_{n_L}^T \mbf{R}(t) }{r_{n\low{L}}\left(t\right)} \! & \! \mbf{0} \! & \! 0 \! & \! \frac{1}{r_{n_L}\left(t\right)} \! & \! 0 \! & \! 0 \! & \! 0\\
			\mbf{a}^T\!(t) \mbf{R}^T\!(t) \! & \! \mbf{0} \! & \! \mbf{0} \! & \! 0 \! & \! 0 \! & \! 1 \! & \! 0 \! & \! 0\\
			\mbf{0} \! & \! 2 \mbf{a}^T(t) \! & \! \mbf{0} \! & \! 0 \! & \! 0 \! & \! 0 \! & \! 3 \! & \! 0\\
			\mbf{0} \! & \!\mbf{0} \! & \! \mbf{a}^T(t) \! & \! 0 \! & \! 0 \! & \! 0 \! & \! 0 \! & \! 1\\
			\mbf{0} \! & \! \mbf{0} \! & \! \mbf{0} \! & \! 0 \! & \! 0 \! & \! 0 \! & \! 0 \! & \! 0
		\end{array}\!
	\right]\!\!, &
\end{eqnarray*}

Including system disturbances and sensor noise to tune the Kalman
filter gives the final system dynamics
\[
	\left\{
		\begin{array}{l}
			\dot{\mbg{\chi}}(t) = \mbg{\cal A}(t) \mbg{\chi}(t) + \mbf{B}(t) \mbf{a}(t) + \mbf{n}_x(t)\\
			\mbf{y}(t) = \mbf{C}(t) \mbg{\chi}(t) + \mbf{n}_y(t)
		\end{array}
	\right.,
\]
where it is assumed that $\mbf{n}_x$ and $\mbf{n}_y$ are
uncorrelated zero-mean Gaussian noise, with
${\rm{E}}\left[ \mbf{n}_x(t) \mbf{n}_x^T\left(\tau\right) \right] = \mbf{Q}_x \delta\left(t-\tau\right)$
and
${\rm{E}}\left[ \mbf{n}_y\left(t\right) \mbf{n}_y^T \left(\tau\right) \right] = \mbf{Q}_y \delta\left(t-\tau\right).$

\subsection{Practical Considerations}
\label{sec:FD:PC}

It was assumed in the paper that there were at least 4 noncoplanar
landmarks, in order to fit the configuration of a Long Baseline
acoustic setup. Nevertheless, the proposed solution is general and
can be applied to any number of landmarks, including the case of
single range measurements. This last setup was already studied by
the authors, see \cite{paper:Batista:ECC:2009} and
\cite{paper:Batista:ACC:2010}. The remaining cases of 2 and 3
landmarks differ only in the observability analysis. When there is
a single range measurement, it was shown in
\cite{paper:Batista:ACC:2010} that the corresponding LTV system would
be observable if and only if the set of functions

{\small
\begin{eqnarray*}
	& {\cal F} = \left\{ t-t_0,\, \left(t-t_0\right)^2,\, \left(t-t_0\right)^3,\, \left(t-t_0\right)^4, \right.&
	\nonumber\\
	& \left. p_1(t) - p_1\left(t_0\right),\, p_2(t) - p_2\left(t_0\right),\, p_3(t) - p_3\left(t_0\right), \right. &
	\nonumber\\
	& \left. \left(t-t_0\right) \left[ p_1(t) - p_1\left(t_0\right) \right],\, \left(t-t_0\right) \left[ p_2(t) - p_2\left(t_0\right) \right], \right. &
	\nonumber\\
	& \left. \left(t-t_0\right) \left[ p_3(t) - p_3\left(t_0\right) \right], \, \left(t-t_0\right)^2 \left[ p_1(t) - p_1\left(t_0\right) \right], \right. &
	\nonumber\\
	& \left. \left(t-t_0\right)^2 \left[ p_2(t) - p_2\left(t_0\right) \right],\, \left(t-t_0\right)^2 \left[ p_3(t) - p_3\left(t_0\right) \right] \right\}&
\end{eqnarray*}}\noindent
is linearly independent on $\left[t_0,t_f\right]$. The cases of
2 and 3 landmarks require less demanding conditions so that the
corresponding LTV systems are observable. These conditions will be
presented elsewhere. Nevertheless, the design follows the same
steps.

\section{Simulation Results}
\label{sec:SR}

In order to evaluate the performance achieved with the proposed
navigation solution, simulations were carried out using a kinematic
model for an underwater vehicle. The fact that the full nonlinear
dynamics of the vehicle are not considered is not a drawback as the
proposed filter relies solely on the vehicle kinematics, which are
exact. Therefore, the proposed solution applies to any underwater
vehicle, independently of the particular dynamics.

\begin{figure}[htbp]
	\center
	{
		\includegraphics[width=0.45\textwidth, keepaspectratio]{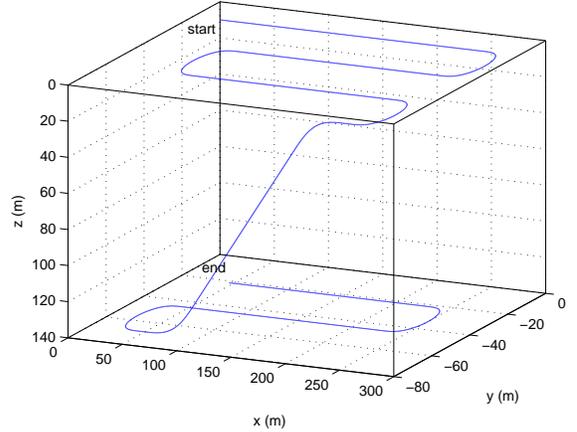}
	}
	\caption{Trajectory described by the vehicle}
	\label{fig:SR:Trajectory}
\end{figure}

The trajectory described by the vehicle is shown in Fig.
\ref{fig:SR:Trajectory}.

Sensor noise was considered for all sensors. In particular, the
range, acceleration, and angular velocity measurements are assumed
to be corrupted by additive uncorrelated zero-mean white Gaussian
noise, with standard deviations of $1 \, \rm{m}$,
$2 \times 10^{-3} \, \rm{m/s^2}$, and $0.05$ \textdegree/s,
respectively. The attitude, parameterized by roll, pitch, and yaw
Euler angles, was also assumed to be corrupted by zero-mean additive
white Gaussian noise, with standard deviation of $0.03${\textdegree}
for the roll and pitch and $0.3${\textdegree} for the yaw.

The LBL configuration is composed of 4 acoustic transponders and
their positions are
\[
	\mbf{s}_1
	=
	\left[
		\begin{array}{c}
			0\\
			0\\
			1000
		\end{array}
	\right]
	\rm{(m)},
	\;
	\mbf{s}_2
	=
	\left[
		\begin{array}{c}
			1000\\
			0\\
			1000
		\end{array}
	\right]
	\rm{(m)},
\]
\[
	\mbf{s}_3
	=
	\left[
		\begin{array}{c}
			0\\
			1000\\
			1000
		\end{array}
	\right]
	\rm{(m)},
	\;
	\mbf{s}_4
	=
	\left[
		\begin{array}{c}
			0\\
			0\\
			500
		\end{array}
	\right]
	\rm{(m)},
\]
which satisfy Assumption \ref{assumption:PS:LBL Configuration}.

To tune the Kalman filter, the state disturbance intensity matrix
was chosen as $\mbf{Q}_x = 10^{-5} \mbf{I}$
and the output noise intensity matrix as
$\mbf{Q}_y = \rm{diag} \left(1,1,1,1,2,2,2,2,2,2\right)$.
The initial conditions were set to zero for the position and
velocity. The acceleration of gravity was initialized closed to the
true value, with $\left[0\;0\;10\right]^T \rm{m/s^2}$ as the
attitude is measured and the magnitude of the acceleration of
gravity is usually known. Notice that it would be possible to
initialize the position with a close estimate obtained from the
inversion of the first set of LBL range measurements. The states
corresponding to the range measurements were initialized with the
first set of measurements while the remaining states were set to
zero, apart from $x_{11}$, which corresponds to the square of the
magnitude of the acceleration of gravity, which was initialized with
100.

\begin{figure}[htbp]
	\center
	{
		\includegraphics[width=0.45\textwidth, keepaspectratio]{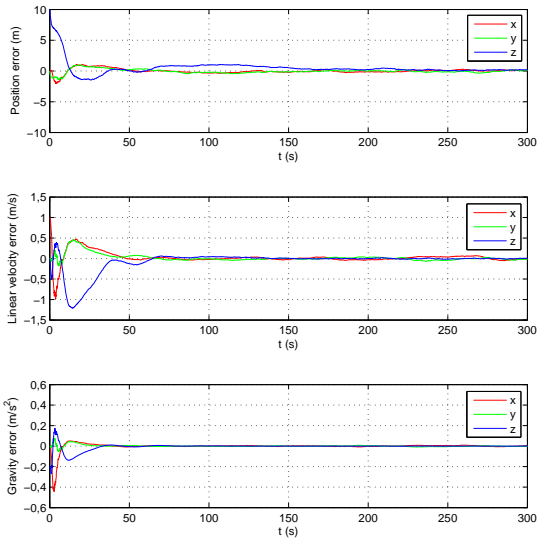}
	}
	\caption{Initial convergence of the position, velocity, and acceleration of gravity error}
	\label{fig:SR:PS:Error variables - Initial Conv PVG}
\end{figure}

\begin{figure}[htbp]
	\center
	{
		\includegraphics[width=0.45\textwidth, keepaspectratio]{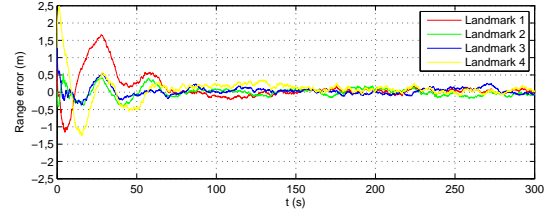}
	}
	\caption{Initial convergence of the range errors}
	\label{fig:SR:PS:Error variables - Initial Conv Range}
\end{figure}

\begin{figure}[htbp]
	\center
	{
		\includegraphics[width=0.45\textwidth, keepaspectratio]{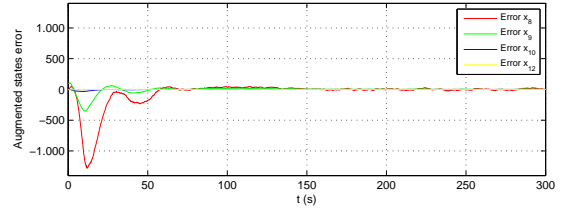}
	}
	\caption{Initial convergence of the error of the augmented states $x_8$, $x_9$, $x_{10}$, and $x_{11}$}
	\label{fig:SR:PS:Error variables - Initial Conv AS}
\end{figure}

The initial evolution of the position, velocity, and acceleration of
gravity errors is depicted in Fig.
\ref{fig:SR:PS:Error variables - Initial Conv PVG},
whereas the initial evolution of the range errors is shown in Fig.
\ref{fig:SR:PS:Error variables - Initial Conv Range}. The initial
convergence of the remaining state errors is shown
in Fig. \ref{fig:SR:PS:Error variables - Initial Conv AS}. As it can
be seen from the various plots, the convergence rate of the filter
is quite fast.

\begin{figure}[htbp]
	\center
	{
		\includegraphics[width=0.45\textwidth, keepaspectratio]{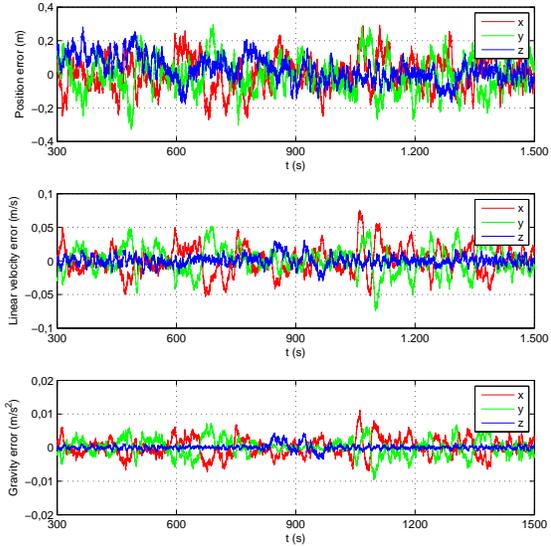}
	}
	\caption{Detailed evolution of the position, velocity, and acceleration of gravity error}
	\label{fig:SR:PS:Error variables - Detail PVG}
\end{figure}

In order to better illustrate the performance achieved with the
proposed solution, the steady-state errors of the position,
velocity, and acceleration of gravity are shown in Fig.
\ref{fig:SR:PS:Error variables - Detail PVG}. Notice that the errors
are confined to very tight intervals, in spite of the realistic
measurement noise. 

The proposed filter was compared against two different solutions:
the first is the well-known Extended Kalman Filter;
the second consists in applying the linear Kalman filter proposed in
\cite{paper:Batista:ACC:2009:1} using a position algebraic
estimate obtained from the range measurements to feed the filter.
Unfortunately, it is not possible to show the results here due to
the lack of space, but in short, the performance of the proposed
solution is similar to the EKF and outperforms the one based on the
algebraic inversion of the problem. These results will be shown
in an extended version of the paper.

\section{Conclusions}
\label{sec:C}

This paper presented a novel Long Baseline (LBL) position and
velocity navigation filter for underwater vehicles based directly on
the sensor measurements. Traditional solutions resort either to the
Extended Kalman Filter (EKF) or to solutions based on position
algebraic estimates obtained from the range measurements. The
solution presented in the paper departs from previous approaches
as the range measurements are explicitly embedded in the filter
design, therefore avoiding inversion algorithms. Moreover, the
nonlinear system dynamics are considered to their full extent and no
linearizations are carried out whatsoever, which allows to show that
the filter error dynamics are globally asymptotically stable.
State augmentation is at the core of the proposed framework. Indeed,
the proposed filter has $13 + n\low{L}$ states, which compares to
a minimum of 9 states for traditional solutions. It is the opinion
of the authors that the achieved performance coupled with the
guarantee of global asymptotic stability justify the additional
computational complexity, at least for those missions where extreme
constraints on power do not exist. 

Under simulation environment it was shown that the filter achieves
similar performance to the Extended Kalman Filter (EKF) and
outperforms linear position and velocity filters based on position
algebraic estimates obtained directly from the range measurements.
The advantage over the EKF is however clear due to the GAS
property of the proposed solution, which is not guaranteed for the
EKF.

Finally, although sensor bias, in particular, accelerometer bias, is
a very important problem, in this paper it is assumed, for the sake
of clarity of presentation and simplicity of design, that the
accelerometer has been previously calibrated. Nevertheless, it is
easy to include the accelerometer bias in the system design, which
leads to observability conditions so that it is possible to estimate
both the gravity and the bias, according to previous work by the
authors. The overall system design and observability analysis will
be presented elsewhere.

\bibliography{References}

\end{document}